\documentclass[11pt]{article}
\usepackage{a4}
\usepackage{amsthm, amsmath, amsfonts, enumitem}
\usepackage{graphicx}
\usepackage{caption}
\usepackage{titlesec}
\graphicspath{ {./Paperfigures/} }
\usepackage{epsfig}
\usepackage{caption}  
\usepackage{capt-of}
\usepackage{url}
\usepackage{tikz}
\usepackage{enumitem}
\usepackage{multicol}
\usepackage{thmtools}
\usepackage{thm-restate}
\usepackage{url}
\usepackage{cleveref}
\newtheorem{theorem}{Theorem}

\newtheorem{lemma}[theorem]{Lemma}

\title{Locally Hamiltonian graphs and minimal size of maximal graphs on a surface}
\author{James Davies\thanks{Department of Combinatorics and Optimization, University of Waterloo, Waterloo, Canada. E-mail: \texttt{jgdavies@uwaterloo.ca}.} \ and Carsten Thomassen\thanks{Department of Applied Mathematics and Computer Science, Technical University of Denmark, DK-2800 Lyngby, Denmark. Email: \texttt{ctho@dtu.dk}.
This work was done while the second author held the Dean's Distinguished Visiting Professorship at the University of Waterloo Fall 2019.}}
\date{}

\begin{document}

\maketitle

\begin{abstract}
	We prove that every locally Hamiltonian graph with $n\ge 3$ vertices and possibly with multiple edges has at least $3n-6$ edges with equality if and only if it triangulates the sphere. As a consequence, every edge-maximal embedding of a graph $G$ graph on some 2-dimensional surface $\Sigma$ (not necessarily compact) has at least $3n-6$ edges with equality if and only if $G$ also triangulates the sphere.
	If, in addition, $G$ is simple, then for each vertex $v$, the cyclic ordering of the edges around $v$ on $\Sigma$ is the same as the clockwise or anti-clockwise orientation around $v$ on the sphere. If $G$ contains no complete graph on 4 vertices and has at least 4 vertices, then the face-boundaries are the same in the two embeddings.
\end{abstract}

\section{Introduction}

In 1974 Kainen \cite{kainen1974some} posed as the first open problem in his survey the following question: By how many edges can an edge-maximal graph embeddable in a surface of Euler genus $g$ be short of a triangulation? Unlike planar graphs, a non-complete edge-maximal graph on a given surface does not necessarily triangulate the surface \cite{davies2019,franklin1934six,harary1974maximal,huneke1978minimum,jungerman1980minimal,ringel1955man}. The first author and Pfender \cite{davies2019} proved that for surfaces of Euler genus at most 2, there are exactly two such examples, being $K_7-e$ on the Klein bottle and $K_8-E(C_5)$ on the torus. In contrast to this they further showed that for orientable surfaces of genus $g\ge 2$, there exist infinitely many such graphs that are $\lfloor \frac{g}{2} \rfloor$ edges short of a triangulation.

McDiarmid and Wood \cite{mcdiarmid2018edge} pointed out that there is a related, equally natural question: By how many edges can an edge-maximal graph embedded in a surface of Euler genus $g$ be short of a triangulation? Indeed if a graph is edge-maximal embeddable in a surface, then every embedding is edge-maximal. They gave a first answer to both problems with the upper bound of $84g$.

In this paper we provide the best possible bound to McDiarmid and Wood's question and characterize the extremal graphs. It is perhaps worth noting that the answer does not depend on the genus. It does not even depend on the genus being finite. Our result also holds for graphs with multiple edges where “edge-maximal” means that we cannot add an edge between non-neighbours. We prove that an edge-maximal embedding of a $n\ge 3$ vertex graph $G$ on a surface $\Sigma$ has at least $3n-6$ edges, with equality if and only if $G$ has an embedding triangulating the sphere. Furthermore, when equality holds and $G$ is simple, then for each vertex $v$, the cyclic ordering of the edges around $v$ on $\Sigma$ is the same as the clockwise or anti-clockwise orientation around $v$ on the sphere. If $G$ contains no complete graph on 4 vertices and has at least 4 vertices, then the face-boundaries are the same in the two embeddings.
This can be seen as a generalisation of the classical folklore result that edge-maximal planar embeddings triangulate the plane. As a consequence this also provides an upper bound of $3g-1$ to Kainen's question for each $g>0$.

For simple graphs, the result follows from an earlier result by Skupie{\'n} \cite{skupien1966locally}: Every connected, locally Hamiltonian simple graph on $n\ge 3$ vertices has at least $3n-6$ edges, with equality if and only if the graph triangulates the sphere\footnote{Skupie{\'n} points out in \cite{law2008locally} that the proof in \cite{skupien1966locally} is incomplete and states the theorem as an open problem. We point out that marginal additions to \cite{skupien1966locally} completes the proof of the stated theorem.
}. 
By adapting Skupie{\'n}'s proof, we prove a generalisation for graphs with multiple edges.

It is an easy consequence of Kuratowski’s theorem that a 3-connected graph distinct from $K_5$ is planar if and only if it contains no subdivision of $K_{3,3}$. The second author \cite{thomassen1981kuratowski} conjectured that a 4-connected simple graph with at least $3n-6$ edges is planar if and only if it contains no subdivision of $K_5$. This was motivated by, and would imply, the conjecture of Dirac \cite{dirac1964homomorphism} that every simple graph with $n$ vertices and more than $3n-6$ edges contains a subdivision of $K_5$. Both of these conjectures were proved by Mader \cite{mader19983n,mader2005graphs}. Combining Mader's theorem with the theorem in the present paper we get the following: If a 4-connected simple graph $G$ has an edge-maximal embedding on some surface $\Sigma$, then it contains a subdivision of $K_5$, unless each face of $G$ on $\Sigma$ is bounded by a triangle, and the replacement of each face by a disc results in the sphere triangulated by $G$.

\section{Preliminaries}

The notation is essentially that of \cite{mohar2001graphs}. A graph has no loops except in the last section, but may have multiple edges. A vertex is \emph{simple} if it is not incident with any multiple edges. A graph is \emph{simple} if all its vertices are simple. A \emph{surface} is an arcwise connected Hausdorff space which is locally homeomorphic to an open disc in the plane. We shall not assume that $\Sigma$ is compact. If $G$ is an abstract graph, then an \emph{embedding} of $G$ on $\Sigma$ is a graph $G^*$ isomorphic to $G$ such that each vertex of $G^*$ is a point in $S$ and each edge in $G^*$ is a simple arc joining its two ends such that two edges may only intersect in a common end. The edges leaving a vertex $v$ have two cyclic orderings (one being the reverse of the other). We choose one of these and call it the \emph{clockwise orientation} around $v$ in $G^*$, while the other we call the \emph{anti-clockwise orientation}.

An embedding $G^*$ on a surface $\Sigma$ is \emph{edge-maximal} if, for any two non-adjacent vertices $x,y$, it is not possible to add an edge $xy$. More precisely, $\Sigma$ does not contain a simple arc from $x$ to $y$ that does not contain any other point of $G^*$.

If $G^*$ is edge-maximal, it clearly satisfies the following.

\begin{itemize}
	\item The graph $G^*$ is connected. If $G^*$ has at least 3 vertices and $v$ is a vertex of $G^*$, then $v$ has at least two neighbours.
	
	For otherwise, we can add an edge close to an edge from $v$ to a neighbour $u$ and then close to an edge from $u$ to a neighbour $w$, thereby adding a new edge $vw$, contradicting the maximality.
	
	\item If $vv_1,vv_2,…,vv_d$ are the edges incident with $v$ in clockwise order (where the indices are expressed modulo $d$), then $G^*$ contains the edge $v_iv_{i+1}$ whenever the vertices $v_i,v_{i+1}$ are distinct.
	
	For otherwise, we could add that edge close to the path $v_ivv_{i+1}$.
\end{itemize}

Motivated by this we define a vertex $v$ of a graph to be \emph{locally Hamiltonian} if the edges incident with $v$ have a cyclic ordering (called a \emph{Hamiltonian ordering}) $vv_1,vv_2,…,vv_d$ so that $G$ contains the edge $v_iv_{i+1}$ whenever the vertices $v_i,v_{i+1}$ are distinct. A graph is \emph{locally Hamiltonian} if every vertex is locally Hamiltonian. Observe that if an embedding $G^*$ is edge-maximal, then the graph $G$ is locally Hamiltonian. Note that, if we replace an edge in a locally Hamiltonian graph by a multiple edge, the resulting graph remains locally Hamiltonian. However, the converse is not true as shown by the planar triangulations with multiple edges defined below.

A \emph{planar triangulation} or a \emph{triangulation of the sphere} is a graph that can be embedded in the plane or the sphere such that every face is bounded by a 3-cycle. Thus $K_3$ is the only planar triangulation on 3 vertices, and $K_4$ and the graph consisting of a 2-cycle and two simple vertices of degree 2 joined to both vertices of the 2-cycle are the only planar triangulations on 4 vertices.

We shall use a lemma on planar triangulations that may have multiple edges.

\begin{lemma}
	Let $T$ be a planar triangulation that contains a cycle $C$ that either has length 2 or is a non-facial (that is, separating) 3-cycle. If all vertices in the interior (respectively exterior) have degree at least 4, then there exist at least 2 vertices $a,b$ in the interior (respectively exterior), such that each is simple, has degree at most 5, has at most two non-simple neighbours, and such that none is contained in a non-facial 3-cycle.
\end{lemma}

\begin{proof}
	Assume without loss of generality that all vertices in the interior have degree at least 4. As $T$ is a triangulation, if $C$ is a 2-cycle then it must be non-facial and separating. We may further assume that the interior of $C$ does not contain another such cycle $C^*$. It is enough to show that there exist vertices $a$ and $b$ in the interior of $C$ satisfying the conclusion of Lemma 1.
	
	Now the vertices in the interior are all simple as there are no 2-cycles in the interior of $C$. Each vertex of $C$ is adjacent to at least one additional vertex in the interior of $C$ and every vertex in the interior has degree at least 4. In particular there must be at least 2 vertices in the interior. As there is no separating 3-cycle, every vertex in the interior is adjacent to at most 2 vertices of $C$ and so no vertex has more than 2 non-simple neighbours. Now by Euler's formula there exist at least 2 vertices in the interior of $C$ with degree at most 5. So we may choose 2 of these to be $a$ and $b$ as required.
\end{proof}

\section{Locally Hamiltonian graphs}

In this section we generalise Skupie{\'n}'s theorem to graphs with multiple edges.

\begin{theorem}
	Let $G$ be a connected, locally Hamiltonian graph with $n\ge 3$ vertices and at most $3n-6$ edges. Then $G$ has precisely $3n-6$ edges and is isomorphic to a triangulation of the sphere.
\end{theorem}

\begin{proof}
	We adapt Skupie{\'n}'s proof. The proof is by induction. If $n=3$, $G$ contains a 3-cycle, so assume $n>3$. Note that $G$ contains a vertex of degree at most 5 as the number of edges is less than $3n$.
	
	Assume now (reduction ad absurdum) that Theorem 2 is false. Let $G$ be a counterexample and let $v$ be a vertex of degree $d\le 5$ such that;
	
	\begin{enumerate}[label=(\roman*)]
		\item $n$ is smallest possible and, subject to (i),
		
		\item if possible $d\le 3$ and subject to (i), (ii),
		
		\item the number of edges which are incident with $v$ and part of multiple edges is minimum, and, subject to (i), (ii), (iii),
		
		\item the number of edges joining non-consecutive vertices in $N(v)$ is smallest possible and, subject to (i), (ii), (iii), (iv),
		
		\item the number of non-simple vertices in $N(v)$ is minimum, and subject to (i), (ii), (iii), (iv), (v),
		
		\item the degree $d$ of $v$ is smallest possible.
	\end{enumerate}
	
	Let $v_1,v_2,…,v_d$ be the edges incident with $v$ in cyclic order (where the indices are expressed modulo $d$) such that $G$ contains the edge $v_iv_{i+1}$ for $i=1,2,…,d$ whenever $v_i,v_{i+1}$ are distinct.
	
	Consider first the case $d=2$. As $n>3$ and $G$ is connected and contains the 3-cycle $vv_1v_2v$, one of $v_1,v_2$, say $v_1$, has degree at least 3. Let $e_1$, respectively $e_2$, be the edge preceding, respectively succeeding, $v_1v$ in the cyclic ordering around $v_1$ given by a Hamiltonian ordering in $G[N(v_1)]$. As $v_1$ has degree at least 3, $e_1$ and $e_2$ are distinct. Then $e_1,e_2$ both join $v_1$ to some neighbour of $v$, and that neighbour must be $v_2$. It is easy to see that $G-e_1-v$ is locally Hamiltonian. By the induction hypothesis, $G-e_1-v$ triangulates the sphere. It is easy to extend that triangulation by adding the edges $e_1,vv_1,vv_2$ to give an embedding of $G$ triangulating the sphere.
	
	Consider next the case $d=3$. It is easy to see that $G-v$ is locally Hamiltonian. If $v_2=v_3$, say, then the proof of the case $d=2$ shows that G has at least two edges $e_1,e_2$ between $v_1$ and $v_2$. Repeating the proof of the case $d=2$, we conclude that  $G-e_1-v$ is locally Hamiltonian. By the induction hypothesis, $G$ has $3(n-1)-6 + 4=3n-5>3n-6$ edges,  a contradiction. So we may assume that $v_1,v_2,v_3$ are distinct.
	
	Now $G-v$ is locally Hamiltonian, so by the minimality of $G$, the graph $G-v$ has an embedding $(G-v)^{**}$ triangulating the sphere. We claim that the 3-cycle $v_1v_2v_3v_1$ must be facial since otherwise $G-v_1-v_2-v_3$ would have three components with $v_1$ being adjacent to a vertex of each component, but then $v_1$ wouldn't be locally Hamiltonian. So $v_1v_2v_3v_1$ is facial, and hence $G$ is a triangulation.
	
	We may now assume by (ii), that $G$ has minimum degree at least 4.
	
	Consider next the case $d=4$. If any two of $N(v)$ are neighbours, then $G-v$ is locally Hamiltonian. But then $G-v$ would have only $3(n-1)-7$ edges, contradicting assumption (i). So we may assume that $v_1,v_3$ are distinct and non-adjacent. Hence $v_1,v_2,v_3$ are distinct, and $v_1v_4,v_3$ are distinct (but possibly $v_2=v_4$), and $H=G-v+v_1v_3$ is locally Hamiltonian with an embedding $H^{**}$ triangulating the sphere.
	
	Now if $v_1v_3v_2v_1$ and $v_1v_3v_4v_1$ are both facial in $H^{**}$, then it is straightforward to obtain the triangular embedding of $G$. So we may assume that one is not facial and apply Lemma 1 to deduce by (iii), (iv) and (v) that $v_2$ and $v_4$ are distinct but not adjacent and without loss of generality that $v_1$ is simple. Every Hamiltonian cycle of $N_G(v_1)$ must contain the path $v_2vv_4$, so there is a Hamiltonian cycle of $N_H(v_1)$ containing the path $v_2v_3v_4$. But then both $v_1v_3v_2v_1$ and $v_1v_3v_4v_1$ must be facial in $H^{**}$, a contradiction.
	
	Finally consider the case $d=5$. We add to $N(v)$ say $m$ edges so that any two vertices in $N(v)$ are adjacent and remove the vertex $v$ to obtain a locally Hamiltonian graph $H$. If $m\le 1$, then $H$ has at most $3(n-1) -7$ edges, a contradiction.
	
	Suppose that $m = 2$, then $H$ has at most $3(n-1)-6$ edges and so is there is an embedding $H^{**}$ triangulating the sphere. This triangulation $H^{**}$ must be simple since otherwise Lemma 1 gives a contradiction. If $v$ is non-simple then as $m=2$, we must have that $N(v)=\{u_1,u_2,u_3,u_4\}$ and we may assume that $u_1,u_3$ are non-adjacent in $G$, and hence in $H$ there is just a single edge between $u_1$ and $u_3$. Suppose that $u_1u_3u_2u_1$ and $u_1u_3u_4u_1$ are both facial, then it is easy to obtain an embedding of $G$ with one just two non-triangular faces, a 2-face $vu_iv$ for some $i\in\{1,2,3,4\}$, and a 4-face $a_1a_2a_3a_4a_1$. We may assume without loss of generality that $a_1$ and $a_2$ are simple (because $H^{**}$ is simple and hence $v$ and $u_i$ are the only non-simple vertices in $G$). As the neighbourhood of a simple vertex of a planar graph has at most one Hamiltonian cycle, it follows that $a_1$ is adjacent to $a_3$ and $a_2$ is adjacent to $a_4$. But then we may obtain a planar embedding of $K_5$, a contradiction. Hence one of these two 3-cycles is non-facial and we may apply Lemma 1, to contradict (iii). So $v$ must be simple in $G$. But then $H$ contains a $K_5$, contradicting the fact that $H$ is planar.
	
	So $m\ge 3$ and as a consequence $v$ must be simple. In particular $G[N(v)]$ contains at most two edges besides the 5-cycle $v_1v_2v_3v_4v_5v_1$. So we may choose the notation such that either $G[N(v)]$ contains no edge or the edge $v_2v_4$ or the two edges $v_2v_4,v_3v_5$ or the edges $v_2v_4,v_2v_5$. In the latter case we may assume  by symmetry that some Hamiltonian cycle in $G[N(v_2)]$ does not contain the path $v_3vv_5$.
	
	Let $Q$ be obtained by deleting $v$ and adding the edges $v_1v_3,v_1v_4$. The resulting graph is locally Hamiltonian and hence isomorphic to a triangulation $Q^{**}$ of the sphere. If the 3-cycles $v_1v_2v_3v_1$, $v_1v_3v_4v_1$ and $v_1v_4v_5v_1$ are all facial then $Q^{**}$ can easily be modified to show that $G$ is also a spherical triangulation. So we may assume that at least one is not facial. Now $v_2$ is not adjacent to $v_4$, else we may apply Lemma 1 to contradict (iv). Hence we may assume that the only edges of $G[N(v)]$ are those of the 5-cycle $v_1v_2v_3v_4v_5v_1$. Some vertex of $N(v)$ must be simple, as otherwise there would be a 2-cycle in $Q^{**}$ whose interior contains no vertex of $N(v)$, allowing for a contradiction by Lemma 1 and (v). So without loss of generality we may assume that $v_1$ is simple. Now as $v_2vv_5$ is a path in every Hamiltonian cycle of $N_G(v_1)$, there is a Hamiltonian cycle of $N_Q(v_1)$ containing the path $v_2v_3v_4v_5$. As $v$ is simple, $N_Q(v_1)$ has only one Hamiltonian cycle and therefore $v_1v_2v_3v_1$, $v_1v_3v_4v_1$ and $v_1v_4v_5v_1$ are all facial, a contradiction.
\end{proof}

\section{Maximal embedded graphs with minimum number of edges}

We are now ready to prove our main result on edge-maximal embeddings.

\begin{theorem}
	Let $G$ be a graph with $n\ge 3$ vertices and at most $3n-6$ edges, and let $G^*$ be an embedding of $G$ on a surface $\Sigma$. If $G^*$ is edge-maximal on $\Sigma$, then $G$ is planar and has precisely $3n-6$ edges. Moreover, there is an embedding $G^{**}$ on the sphere such that $G^{**}$ is a triangulation.
	
	If $G$ is simple then, for every vertex $v$ in $G$, the clockwise orientation around $v$ in $G^*$ is either the clockwise or anti-clockwise orientation around $v$ in $G^{**}$.
	If $G$ is simple, $n\ge 4$, and contains no $K_4$, then $G^*$ is obtained from $G^{**}$ by replacing each face by a surface.
\end{theorem}

\begin{proof}
Since $G^*$ is edge-maximal on $\Sigma$, it follows that $G$ is locally Hamiltonian. By Theorem 1, $G$ has at least $3n-6$ edges, and if equality holds, $G$ has an embedding $G^{**}$ on the sphere such that $G^{**}$ is a triangulation.

If $G$ is simple, the embedding of $G^{**}$ is unique except that we can interchange “clockwise” with “anti-clockwise” for every vertex. The clockwise orientation in $G^*$ is given by some Hamiltonian cycle in $G[N(v)]$. But $G^{**}[N(v)]$ has precisely one Hamiltonian cycle. So, the clockwise orientation at $v$ in $G^*$ equals the clockwise or anticlockwise orientation at $v$ in $G^{**}$.

Suppose now that $G$ is simple, $n  \geq 4$, and contains no $K_4$. Note that $G$ has minimum degree at least 3. Each edge of $G$ appears on facial walks of $G^*$  exactly twice. Since $G^*$ is edge-maximal, each face boundary induces a complete graph. In particular, each facial walk induces a complete graph. Since $G^*$ is simple and contains no $K_4$, it follows that each
facial walk in $G^*$  has exactly 3 vertices and then also exactly 3 edges. Consider now a vertex $v$ of degree $d\ge 3$, say. In $G^{**}$ there is a unique Hamiltonian cycle $v_1v_2…v_dv_1$ in $G[N(v)]$, defined by the clockwise ordering in $G^{**}$ of the edges incident with $v$. As this is also the clockwise or anticlockwise ordering in $G^*$ of the edges incident with $v$, all $d$ facial walks in $G^*$ containing $v$ are precisely the facial walks in $G^{**}$ containing $v$. So all facial walks in $G^{**}$ are facial walks in $G^*$. Hence the facial walks in $G^*$ are precisely the facial triangles in $G^{**}$. In $G^{**}$ a face has only one facial walk. So, there only remains the question if a face in $G^*$ may have more than one facial walk. But, this is impossible because every face boundary is complete and has therefore only 3 vertices. Hence $G^*$  is obtained from $G^{**}$  by replacing each face by a surface.
\end{proof}

Note that the orientation around vertices in graphs with multiple edges is not as simple. Indeed, the planar triangulation with vertices $a,b,c,u,v$ and edges $uv,ua,ub,uv,uc, va,vb,vc,ab$ has an edge-maximal embedding in the projective plane such that the clockwise orientation around $u$ is $uv,ub,ua,uv,uc$, that is the same as on the sphere except that $ua,ub$ have been replaced by $ub,ua$. Similarly $K_4$ has an embedding in the projective plane with three 4-faces, and so this embedding can not be obtained by replacing faces of a spherical embedding by other surfaces.

\section{Maximal embedded graphs with loops}

In this section we allow loops. Remarkably, edge-maximal embeddings of simple graphs and of graphs with multiple edges behave in much the same way. It is natural to ask if graphs with loops could also behave similarly. However, the answer for edge-maximal embeddings of graph with loops is rather different.

Let $K_n^{o}$ denote the complete graph with a single additional loop on some vertex, we call this vertex the \emph{loop vertex}. We call the two graphs $K_2^o$ and $K_3^o$ \emph{petals}. A graph $G$ is a \emph{flower} if there exists a sequence of graphs $G_0,G_1,\dots ,G_n$ such that, $G_0\in \{K_2,K_3\}$, and for each $i\in \{1,\dots ,n\}$, $G_i$ is obtained by identifying a vertex of $G_{i-1}$ with a loop vertex of a petal and $G_n=G$.

Notice that flowers are planar, however they may have faces with arbitrarily many edges. We leave as a simple inductive exercise to the reader to prove that an edge-maximal embedding of a graph $G$ (possibly with loops) on $n\ge 2$ vertices has at least $2n-3$ edges, with equality if and only if $G$ is a flower.

\bibliographystyle{abbrv}
\bibliography{References.bib}

\end{document}